\documentclass[10pt]{article}
\usepackage{repstyle}
\usepackage{latexsym,epsf,amssymb,subeqn}
\usepackage{multirow}
\usepackage[hidelinks]{hyperref}

\newcommand{\eq}{\begin{equation}}
\newcommand{\eeq}{\end{equation}}
\newcommand{\R}{\mathbb R}

\newcommand{\Oh}{{\mathcal O}}
\newcommand{\I}{{\mathcal I}}
\newcommand{\eps}{{\varepsilon}}
\newcommand{\fhi}{\varphi}
\newcommand{\half}{{1\over2}}
\newcommand{\sfrac}[2]{\mbox{\large{$\frac{#1}{#2}$}}}
\newcommand{\mfrac}[2]{\mbox{\Large{$\frac{#1}{#2}$}}}
\newcommand{\sDelta}{{\mbox{\footnotesize{$\Delta$}}}}
\newcommand{\dt}{\sDelta t}
\newcommand{\dx}{\sDelta x}
\newcommand{\dy}{\sDelta y}

\newcommand{\bb}[1]{\mbox{\boldmath $#1$}}
\newcommand{\uZ}{\underline{Z}}

\begin{document}

\title{ \sc
Error Analysis of Explicit Partitioned Runge-Kutta Schemes
for Conservation~Laws
}
\author{\sc
Willem Hundsdorfer%
\thanks{CWI, PO\,Box 94079, 1090-GB Amsterdam,
The Netherlands ({\tt willem.hundsdorfer@cwi.nl}).
The work of this author is supported by
Award No.\ FIC/2010/05 from King Abdullah University of Science and
Technology (KAUST).}
,
David I.\ Ketcheson%
\thanks{
Division Mathematical and Computer Sciences \& Engineering,
King Abdullah University of Science \& Technolog,y
(KAUST), P.O. Box 4700, Thuwal 23955, Saudi Arabia
({\tt david.ketcheson@kaust.edu.sa}). The work of this author
is supported by Award No.\ FIC/2010/05 from King Abdullah University
of Science and Technology (KAUST).}
,
Igor Savostianov%
\thanks{CWI, PO\,Box 94079, 1090-GB Amsterdam,
The Netherlands ({\tt igor.savostianov@gmail.com}).
The work of this author has been supported by
Award No.\ FIC/2010/05 from King Abdullah University of Science and
Technology (KAUST).}
}

\date{}
\maketitle

\begin{abstract} \noindent 
An error analysis is presented for explicit partitioned Runge-Kutta methods 
and multirate methods applied to conservation laws. The interfaces, 
across which different methods or time steps are used, lead to order 
reduction of the schemes.
Along with cell-based decompositions, also flux-based decompositions
are studied. In the latter case mass conservation is guaranteed, but
it will be seen that the accuracy may deteriorate.

\end{abstract}

\section{Introduction}
\label{Sect:Intro}

Spatial discretization of partial differential equations (PDEs)
lead to systems of ordinary differential equations (ODEs), the
so-called semi-discrete systems. 
In this paper we will consider explicit time stepping schemes applied 
to conservation laws $u_t + \nabla\cdot f(u) = 0$ with a given
spatial discretization.  
The CFL stability condition bounds the time step in terms of the
ratio of local (spatial) mesh width and characteristic speeds.
If either of these factors varies substantially, 
it is natural to use local time steps
that match the local convective velocity or spatial mesh width.
Schemes in which different time steps are used over different
parts of the spatial grid are referred to as multirate schemes. 
Such schemes can be studied in the more general setting of partitioned 
or additive Runge-Kutta methods.

Discontinuities arise in the solution to nonlinear conservation laws,
often leading to numerical oscillations or unphysical values.  Thus 
monotonicity properties and maximum principles become important. 
Step-size restrictions for monotonicity for partitioned Runge-Kutta methods
have been studied in \cite{HMS13}.
In these notes we will consider the accuracy of the methods, assuming the
solution to be sufficiently smooth.  For conservation laws this means 
that accuracy is studied away from shocks.

The classical order of a numerical ODE solver is often larger when applied 
to non-stiff ODEs than when applied to PDEs, where one considers 
time step $\dt$ and spatial mesh width $\dx$ tending to zero simultaneously.  
This phenomenon, known as {\it order reduction}, will be analyzed in this 
paper for partitioned Runge-Kutta methods and multirate methods.  

The system of ODEs in $\R^m$, with given initial value, will be written as
\eq
\label{eq:ODE}
u'(t) = F(t,u(t)) \,, \qquad u(0) = u_0 \,.
\eeq
In our applications, this ODE system will be a semi-discrete system
obtained from a conservation law by a finite difference or finite volume 
discretization in space.  Each component $u_j(t)$ of the vector 
$u(t) = [u_j(t)] \in \R^m$ then stands for an approximation at time $t$ 
to the pointwise or average value of the PDE solution at $x_j$,
$j=1,2,\ldots,m$, and $F$ is the spatial discretization operator.

For the time integration of the semi-discrete system we will consider
partitioned methods based on a decomposition of $F$,
\eq
\label{eq:split}
F(t,v) = F_1(t,v) + F_2(t,v) + \cdots + F_r(t,v) \,,
\eeq
where each $F_k: \R\times\R^m\rightarrow\R^m$ corresponds to the spatial
discretization operator in a certain region $\Omega_k$ of the spatial PDE
domain $\Omega = \Omega_1\cup\cdots\cup \Omega_r$.

Let $\I = \I_1\cup\cdots\cup \I_r$ be a partitioning of the index set
$\I = \{1,2,\ldots,m\}$, with $j\in\I_k$ if $x_j\in\Omega_k$. 
To define the schemes we consider corresponding diagonal matrices 
$I = I_1 + \cdots + I_r$, where $I$ is the identity matrix and the 
$I_k$ are diagonal with entries zero or one: the $j$-th diagonal
entry of $I_k$ equal to one iff $j\in\I_k$.
Then $F_k = I_k F$ defines a {\it cell-based} decomposition; the function
$F_k$ contains those components of $F$ that correspond to the spatial
region $\Omega_k$.  
Another possibility is to base the decomposition on fluxes, to ensure
mass conservation; such {\it flux-based} decompositions will be discussed 
later in some detail.
Our main interest is in methods that use different step sizes in each spatial 
domain $\Omega_k$.

The outline of this paper is as follows. Multirate methods are conveniently
analyzed in the broader framework of partitioned and additive Runge-Kutta 
methods, which we review in Section~\ref{Sect:PRKmeths}.
In Section~\ref{sect:MR1}, we present some multirate
methods of order one and two, along with simple numerical tests showing
some of their deficiencies.
General expressions for the local errors, that 
can be used to derive error bounds when both $\dt$ and $\dx$ tend to zero,
are given in Section~\ref{sect:errors}. 
Detailed error bounds are found in Section~\ref{sect:cell-based}
for cell-based decomposition and in Section~\ref{Sect:flux-based} for 
flux-based decomposition.  It will be seen that flux-based decompositions often
lead to a lower order of convergence. 
Some conclusions and final remarks are given in Section~7.

\section{Partitioned Runge-Kutta methods}
\label{Sect:PRKmeths}

For a given decomposition (\ref{eq:split}), we consider partitioned
Runge-Kutta methods, giving approximations $u_n \approx u(t_n)$ at
the time levels $t_n =n\dt$, $n\ge0$.
A step from $t_{n}$ to $t_{n+1}$ with an $s$-stage
method reads
\begin{subequations}
\label{eq:PRK}
\setlength{\arraycolsep}{1mm}
\begin{eqnarray}
v_{n,i}  &  =  &  u_{n} + \dt \sum_{k=1}^r \sum_{j=1}^{s}
a_{ij}^{(k)} F_k(t_n + c_j \dt, v_{n,j})  \,, \qquad
i = 1,\ldots,s\,,
\\
u_{n+1} &  =  &  u_{n} + \dt \sum_{k=1}^r \sum_{j=1}^{s}
b_{j}^{(k)} F_k(t_n + c_j \dt, v_{n,j})  \,.
\end{eqnarray}
\end{subequations}
The internal stage vectors $v_{n,i}$, $i=1,\ldots,s$, give approximations
to $u(t_n + c_i \dt)$ at the intermediate time levels.
For applications to conservation laws we will restrict ourselves to
explicit methods, where $a_{ij}^{(k)} =0$ if $j\ge i$.

For general decompositions $F = F_1+\cdots+F_r$, method (\ref{eq:PRK})
is usually called an additive Runge-Kutta method, and the name partitioned 
Runge-Kutta method is often reserved for the case where the decomposition 
has a partitioned structure ($F_k = I_k F$). However, as noted in 
\cite[p.\,153]{ARS97}, any partitioned method can be written as an additive one
(and vice versa) by modifying the right hand side, so we do not
distinguish these two classes of methods.

In this section we will briefly discuss some basic properties of the
partitioned and additive methods. A more extensive discussion is found 
in \cite{HMS13}.

\subsection*{Internal consistency and conservation}

Let $c_i^{(k)} = \sum_{j=1}^{s} a_{ij}^{(k)}$,
$i = 1,\ldots,s$.  If we have
\eq
\label{eq:intcons}
c_i^{(k)} = c_i^{(l)}  \quad \mbox{ for all $1\le k,l\le r$ and
$1\le i\le s$} \,,
\eeq
then the internal vectors $v_{n,i}$ are consistent approximations
to $u(t_n+c_i\dt)$, and the method is {\it internally
consistent}. As will be seen, this is an important property for the
accuracy of the method when applied to ODEs obtained by
semi-discretization.

If (\ref{eq:intcons}) holds, this gives an obvious choice for the
abscissae $c_i$ in (\ref{eq:PRK}). If (\ref{eq:intcons}) is not satisfied,
then we take $c_i = c_i^{(r)}$, $1\le i\le s$, where it is assumed
that the $r$-th Runge-Kutta method used in (\ref{eq:PRK}) is the most
`refined' one.

Apart from consistency, we will also study {\it conservation} of linear
invariants; for example, mass conservation.
Suppose that $h^T = [h_1,\ldots,h_m]$ is such that
$h^T u(t) = \sum_j h_j u_j(t)$ is a conserved quantity for the ODE system
(\ref{eq:ODE}). This will hold for an arbitrary initial value $u_0$ provided
that
\eq
h^T F(t,v) = 0  \qquad \mbox{for all $t\ge0$, $v \in \R^m$} \,.
\eeq
For the partitioned Runge-Kutta scheme we then have
$$
h^T u_{n+1}
= h^T u_{n} + \dt \sum_{k\neq l} \sum_{j=1}^s
\big(b_{j}^{(k)} - b_{j}^{(l)}\big) h^T F_k(t_n+c_j \dt,v_{n,j}) \,,
$$
for any $1\le l\le r$. Therefore, as noted in \cite{CoSa06}, the
discrete conservation property $h^T u_{n+1} = h^T u_{n}$ will be satisfied
provided that
\eq
\label{eq:conserv}
b_j^{(k)} = b_j^{(l)}
\quad \mbox{ for all $1\le k,l\le r$ and
$1\le j\le s$} \,.
\eeq
If the $h_j$ represent lengths, areas or volumes of cells, this is
often called {\it mass conservation}.
Of course, if $h^T F_k(t,v) \equiv 0$ for all $1\le k\le r$, then
the conservation property will always hold, even if (\ref{eq:conserv})
is not satisfied. This will be valid for decompositions of $F$ that
are based on fluxes.

\subsection*{Order conditions}  
The order conditions for partitioned Runge-Kutta methods applied to
non-stiff problems are found e.g.\ in \cite[Thm.\,I.15.9]{HNW93} for
$r=2$. This order will be denoted by $p$. As we will see, it often
does not correspond to the order of convergence for semi-discrete ODE
systems, and therefore $p$ is usually referred to as the {\it classical order}.

To write the order conditions in a compact way, let the coefficients of
the method be contained in $A_k = [a_{ij}^{(k)}] \in \R^{s\times s}$ and
$b_k = [b_i^{(k)}] \in \R^s$, and set $e = [1,\ldots,1]^T \in \R^s$.
The conditions for order $p$ up to $3$ are
\begin{subequations}
\label{eq:OrderCond}
\setlength{\arraycolsep}{1mm}
\begin{eqnarray}
p = 1:  & &  b_k^T e = 1 \qquad
\mbox{for $\; k = 1,\ldots,r$} \,,
\\
p = 2:  & &  b_k^T A_l \, e = \sfrac{1}{2} \qquad
\mbox{for $\; k, l = 1,\ldots,r$} \,,
\\
p = 3:  & &  b_k^T C_{l_1} A_{l_2} e = \sfrac{1}{3} \,, \quad
b_k^T A_{l_1} A_{l_2} e = \sfrac{1}{6}  \quad
\mbox{for $\; k, l_1, l_2 = 1,\ldots,r$} \,,
\end{eqnarray}
\end{subequations}
where $C_l = \mbox{diag}(A_l e)$.

For semi-discrete ODE systems obtained from a PDE, the accuracy of the
internal stage vectors $v_{n,i} \approx u(t_n + c_i \dt)$ is also of
importance.
The component-wise powers of $c = [c_i] = A_r e$ are denoted by
$c^j = [c_i^{\,j}]$, and $c^0 = e$.
The method is said to have {\it stage order} $q$ if
\eq
\label{eq:StegeOrder}
A_k \, c^j \,=\, \sfrac{1}{j+1} \, c^{j+1} \qquad
\mbox{for $\; j = 0,\ldots,q-1\;$ and $\;k = 1,\ldots,r$} \,.
\eeq
A method is internally consistent if it has stage order $q \ge 1$. 
Furthermore, it is easy to see that an explicit method cannot have $q > 1$. 

Finally, we mention that a necessary condition for having order $p$ is
\eq
b_k^T c^j  \,=\, \sfrac{1}{j+1} \qquad
\mbox{for $\; j = 0,\ldots,p\;$ and $\;k = 1,\ldots,r$} \,.
\eeq


\section{Multirate methods}
\label{sect:MR1}

An important class of methods contained in (\ref{eq:PRK}) are the
multirate methods.
We will consider multirate methods that are based on a single
Runge-Kutta method, such that if $\I_k = \I$ and the
other $\I_l$ are empty, then (\ref{eq:PRK}) reduces to $m_k$
applications of this base method (with step-size $\dt/m_k$,
$m_1 = 1 < m_2 < \cdots < m_r$).
It was shown in \cite{HMS13} that the
conditions for internal consistency (\ref{eq:intcons}) and conservation
of linear invariants (\ref{eq:conserv}) are incompatible for such
multirate schemes.

\subsection{Examples}
\label{sect:Examples}

We consider some explicit multirate schemes that were discussed in
\cite{HMS13}; additional examples can be found e.g.\ in 
\cite{CoSa06,SKAW09,SKAW12,WKG08}.
The schemes in this paper are either based on the forward Euler method
$$
u_{n+1} = u_n + \dt F(t_n,u_n) \,,
$$
or the explicit trapezoidal rule (modified Euler method)
$$
u_{n+1}^* = u_n + \dt F(t_n,u_n) \,, \qquad
u_{n+1} = u_n + \sfrac{1}{2} \dt F(t_n,u_n)
+ \sfrac{1}{2} \dt F(t_{n+1}, u_{n+1}^*) \,.
$$
Furthermore, we take $r=2$, $m_1=1$, $m_2=2$, that is, the local time step
is $\dt$ on $\I_1$ and $\frac{1}{2}\dt$ on $\I_2$.
The coefficients of the schemes are represented by a tableau
$$
\begin{array}{c|c|c}
c  &  A_1  &  A_2
\\
\hline
\rule{0mm}{4mm}
&  b_1^T  &  b_2^T
\end{array}
$$
with $A_k = [a_{ij}^{(k)}] \in \R^{s\times s}$,
$b_k = [b_{i}^{(k)}] \in \R^s$ and $c = [c_i] = A_2 e \in \R^s$.

The scheme with $s = 2$, $p=1$, $q=0$, given by the tableau
\eq
\label{eq:OS1}
\small
\begin{array}{c|cc|cc}
0  &  0  &  0  &  0  &  0    \\
1/2  &  0  &  0  &  1/2 & 0    \\
\hline
 &  1/2 & 1/2  &  1/2 & 1/2
\end{array}
\eeq
is a simple example from Osher \& Sanders \cite{OsSa83}, applied here
with only one level of temporal refinement. We refer to this as the
OS1 scheme.

Another scheme based on forward Euler was given by Tang \& Warnecke
\cite{TaWa06}. It has $s = 2$, $p = q = 1$,
\eq
\label{eq:TW1}
\small
\begin{array}{c|cc|cc}
0  &  0  &  0  &  0  &  0  \\
1/2  &  1/2 & 0  &  1/2 & 0  \\
\hline
 &  1  &  0  &  1/2 & 1/2
\end{array}
\eeq
This scheme is internally consistent but does not conserve linear
invariants because $b_1 \neq b_2$. We will refer to (\ref{eq:TW1}) as
the TW1 scheme.

A second-order scheme of Tang \& Warnecke \cite{TaWa06}, referred to
as the TW2 scheme, is based on the explicit trapezoidal rule.
It has $s = 4$, $p = 2$, $q = 1$,
\eq
\label{eq:TW2}
\small
\begin{array}{c|cccc|cccc}
0    &  0  &  0  &  0  &  0  &  0   &  0  &  0  &  0   \\
1/2  &  1/2 & 0  &  0  &  0  &  1/2 &  0  &  0  &  0   \\
1/2  &  1/4 & 1/4 & 0  &  0  &  1/4 & 1/4 & 0  &  0    \\
1    &  1  &  0  &  0  &  0  &  1/4 & 1/4 & 1/2 & 0    \\
\hline
 &  1/2 & 0  &  0  &  1/2  &  1/4 & 1/4 & 1/4 & 1/4
\end{array}
\eeq

A related scheme, due to Constantinescu \& Sandu \cite{CoSa06}, with
$s = 4 $, $p = 2$, $q = 0$, is given by
\eq
\label{eq:CS2}
\small
\begin{array}{c|cccc|cccc}
0  &  0  &  0  &  0  &  0  &
0  &  0  &  0  &  0    \\
1/2  &  1  &  0  &  0  &  0  &
1/2 & 0  &  0  &  0    \\
1/2  &  0  &  0  &  0  &  0  &
1/4 & 1/4 & 0  &  0    \\
1 &  0  &  0  &  1  &  0  &
1/4 & 1/4 & 1/2 & 0    \\
\hline
 &
1/4 & 1/4 & 1/4 & 1/4  &
1/4 & 1/4 & 1/4 & 1/4
\end{array}
\eeq
This scheme is conservative, but not internally consistent.
We will refer to (\ref{eq:CS2}) as the CS2 scheme.

As a final example we consider the following scheme with $s = 5$,
$p = 2$, $q = 1$,
\eq
\label{eq:SH2}
\small
\begin{array}{c|ccccc|ccccc}
0  &  0  &  0  &  0  &  0  &  0
& 0  &  0  &  0  &  0  &  0     \\
1  &  1  &  0  &  0  &  0  &  0
& 1  &  0  &  0  &  0  &  0     \\
1/2  &  3/8 & 1/8 & 0  &  0  &  0
& 1/2 & 0  &  0  &  0  &  0     \\
1/2  &  3/8 & 1/8 & 0  &  0  &  0
& 1/4 & 0  &  1/4 & 0  &  0     \\
1  &  1/2 & 1/2 & 0  &  0  &  0
& 1/4 & 0  &  1/4 & 1/2 & 0     \\
\hline
 &  1/2 & 1/2 & 0  &  0  &  0
& 1/4 & 0  &  1/4 & 1/4 & 1/4
\end{array}
\eeq
We will refer to this as the SH2 scheme. This scheme has been described in
\cite{HMS13}; it was obtained by adaptation of an implicit (Rosenbrock)
scheme from \cite{SHV07}. Although it looks already a bit complicated,
the idea is simple: first a coarse $\dt$ step is taken with the
explicit trapezoidal rule on the whole index set $\I$ , and then two
refined $\frac{1}{2}\dt$ steps are taken on $\I_2$, using information from the
coarse step by quadratic (Hermite) interpolation at the time level
$t_n + \frac{1}{2} \dt$.

It is important to note that the number of stages $s$ is not a good
measure for the work-load per step.  For example, with the SH2 scheme
we have $s=5$, but neglecting the (small) interface region
only two $F_1$ evaluations and four $F_2$ evaluations are needed per step.

\subsection{Numerical tests in 1D}
\label{Ssect:numtests1}

\subsubsection{Advection with smooth solution}

A convergence analysis of the above multirate schemes,
in the framework of partitioned Runge-Kutta methods,
will be given in the next two sections.
Here we present some simple numerical results for the second-order schemes
that will motivate the analysis. 

To test the accuracy of the schemes we consider the linear
advection equation $u_t + u_x = 0$ on the spatial interval
$\Omega = [0,1]$ with periodic boundary conditions, and time
interval $0<t\le T=1$.
For test purposes a uniform spatial grid is taken, to ensure that interface
effects are not related to the spatial discretization.
The WENO5 finite difference scheme is used; see e.g.\ \cite{Shu09}.  
Further we employ a fixed Courant number $\nu = \dt/\dx = 0.5$, $\dx = 1/m$,
and cell-based splitting $F = I_1 F + I_2 F$, with 
$\I_2 = \{ i : x_i \in 
[\frac{1}{8},\frac{3}{8}]\cup[\frac{5}{8},\frac{7}{8}]\}$.

For this accuracy test a smooth solution $u(x,t) = \sin^2(\pi (x-t))$
is considered. The errors in the maximum norm ($\|v\|_\infty = \max_j|v_j|$)
and discrete $L_1$-norm ($\|v\|_1 = \sum_j \dx_j |v_j|$) are presented in
Table~\ref{Tab:Test1CellBased}.
The entries in the table are the total (absolute) errors
with respect to the exact PDE solution.
In this test the spatial errors are much smaller than the errors
due to time integration with the multirate methods.

\begin{table}[htb]
\vspace{-2mm}
\caption{ \small \label{Tab:Test1CellBased}
Results for the smooth advection problem with the CS2, TW2 and SH2 schemes.
Maximum errors and $L_1$-errors at final time $T = 1$
for various $m$ with fixed Courant number $\dt/\dx = 0.5$, $\dx = 1/m$.
The approximate order of convergence is also given.
}
\small
\begin{center}
\setlength{\tabcolsep}{3mm}
\begin{tabular}{|l||c|c|c|c||c|}
\hline \rule[-1mm]{0mm}{5mm}
\hfill $m$  $\qquad$  & 100    & 200     & 400     & 800 & Order
\\ \hline \hline
\rule{0cm}{4mm}
CS2, $\|\eps_N\|_\infty$
 & $8.22\cdot10^{-4}$ & $2.75\cdot10^{-4}$
 & $1.46\cdot10^{-4}$ & $8.37\cdot10^{-5}$
& 1 \\ \rule{0cm}{3mm}
CS2, $\|\eps_N\|_1$
 & $2.85\cdot10^{-4}$ & $7.81\cdot10^{-5}$
 & $2.09\cdot10^{-5}$ & $5.73\cdot10^{-6}$
& 2 \\[1mm] \hline
\rule{0cm}{4mm}
TW2, $\|\eps_N\|_\infty$
\rule{0cm}{4mm}
 & $3.12\cdot10^{-4}$ & $8.04\cdot10^{-5}$
 & $2.02\cdot10^{-5}$ & $5.05\cdot10^{-6}$
& 2 \\ \rule{0cm}{3mm}
TW2, $\|\eps_N\|_1$
 & $1.98\cdot10^{-4}$ & $5.12\cdot10^{-5}$
 & $1.28\cdot10^{-5}$ & $3.21\cdot10^{-6}$
& 2 \\[1mm] \hline
\rule{0cm}{4mm}
SH2, $\|\eps_N\|_\infty$
 & $3.13\cdot10^{-4}$ & $8.06\cdot10^{-5}$
 & $2.02\cdot10^{-5}$ & $5.05\cdot10^{-6}$
& 2 \\ \rule{0cm}{3mm}
SH2, $\|\eps_N\|_1$
 & $1.99\cdot10^{-4}$ & $5.13\cdot10^{-5}$
 & $1.28\cdot10^{-5}$ & $3.21\cdot10^{-6}$
& 2 \\[1mm] \hline
\end{tabular}
\end{center}
\end{table}

It is seen that with the CS2 scheme we have only first-order convergence
in the maximum norm. The largest errors are localized near the interface 
points; the $L_1$-errors are still second-order.  For the schemes TW2 
and SH2 we have order two convergence also in the maximum norm.

\begin{figure}[htb]
\setlength{\unitlength}{1cm}
\begin{center}
\begin{picture}(9,3.5)
\epsfxsize9cm
\put(0,0){\epsfbox{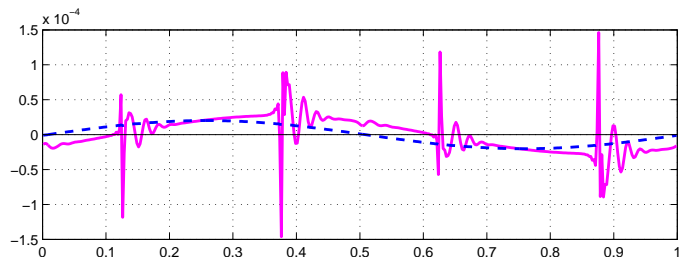}}
\end{picture}
\caption{ \small \label{Fig:A1D}
Errors versus $x_j \in [0,1]$ at final time $T = 1$ for the schemes
CS2 (solid line) and TW2 (dashed line), $m=400$.
}
\vspace{-2mm}
\end{center}
\end{figure}

To see that the largest errors for the CS2 scheme occur indeed
at the interfaces, the errors as function of $x$ at the final time 
$t_n = T = 1$ are displayed in Figure~\ref{Fig:A1D} for $m=400$. 
The (relatively) large errors for the CS2 scheme at the interface points
are clearly visible.  In contrast, the errors for the TW2 scheme show no
visible interface effects; the errors for SH2 were almost the 
same as those for TW2 in this test.

\subsubsection{Shock speeds with Burgers' equation}

The main topic studied in this paper is convergence for smooth solutions. 
Mass conservation will play only a minor role. This conservation property, 
or the lack of it, is of course important for problems with discontinuous 
solutions.

To illustrate this, we apply the CS2, TW2 and SH2 schemes with cell-based 
decomposition to Burgers' equation
\eq
u_t + f(u)_x = 0 \,, \qquad f(u) = \sfrac{1}{2} u^2 \,,
\eeq
with periodic boundary conditions on the spatial region $\Omega = [0,1]$,
and the initial block profile $u(x,0) = 1$ if $x \in [0, \frac{1}{2}]$,
$u(x,0) = 0$ if $x \in [\frac{1}{2}, 1]$.
The shock that starts at $x=\frac{1}{2}$ should be located at $x=\frac{3}{4}$
at the output time $T=\frac{1}{2}$.

A conservative spatial discretization $u_i' = \frac{1}{\Delta x} 
(f_{i-1/2}(u) - f_{i+1/2}(u))$ is used with local Lax-Friedrichs fluxes
\eq
f_{j+\half}(u) \,=\, \sfrac{1}{2} \Big( f(u_{j+\half}^-) + f(u_{j+\half}^+)
+ \alpha_{j+\half} (u_{j+\half}^- - u_{j+\half}^+) \Big) \,,
\eeq
where $\alpha_{j+\half} = \max|f'(v)|$ with $v$ ranging between the states 
$u_{j+1/2}^-$, $u_{j+1/2}^+$ to the left and the right of the cell boundaries,
computed from $u = [u_i]$ with the WENO5 scheme, as in \cite{Shu09}.

The index sets $\I_k$ are changing in time, moving along with the shock.
Let $u_n = [u_i^n]$ denote the numerical solution at time $t_n$; we take 
$\I_1 = \I_1^n = \{i : u_i^n < \frac{1}{8}\}$.
So in the regions with small values of $u_i^n$, corresponding to small
local Courant numbers $\dt |f'(u_i^n)|/\dx$, we use step-size $\dt$;
elsewhere the step-size is $\frac{1}{2}\dt$.

\begin{figure}[htb]
\setlength{\unitlength}{1cm}
\begin{center}
\begin{picture}(13,3.5)
\epsfxsize13cm
\put(0,0){\epsfbox{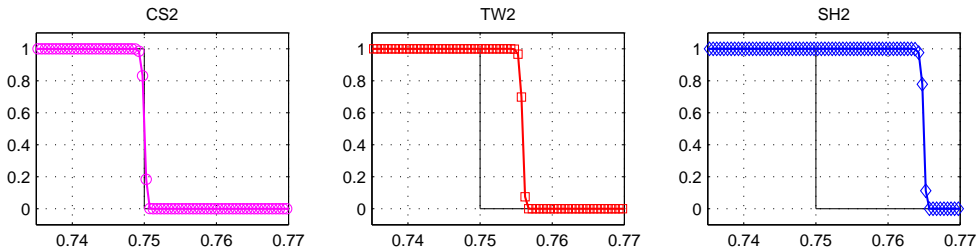}}
\end{picture}
\caption{ \small \label{Fig:BuS}
Shock location for Burgers' equation at time $T=\frac{1}{2}$
for the schemes CS2 (left), TW2 (middle) and SH2 (right), with 
$\I_1 = \I_1^n = \{i : u_i^n < \frac{1}{8}\}$,
$\dt = \dx = 1/m$, $m=2000$.
}
\vspace{-2mm}
\end{center}
\end{figure}

The results with the non-conservative schemes TW2 and SH2 are shown in
Figure~\ref{Fig:BuS}. It is obvious that the lack of conservation leads 
to a shock that moves with a wrong speed; furthermore, the shock does 
not converge to the correct location upon refinement of the grid.  
For the conservative CS2 scheme the shock location is correct.

\section{Local and global discretization errors}
\label{sect:errors}

The local discretization errors of the partitioned methods (\ref{eq:PRK})
will be expressed in terms of derivatives of the functions
\eq
\label{eq:fhi}
\fhi_k(t) \,=\, F_k(t,u(t)) \,.
\eeq
The discretization errors can be studied for nonlinear problems;
see Remark~\ref{Rem:Nonlin}.
However, to avoid unnecessary technical complications we consider only
linear problems with constant coefficients,
\eq
\label{eq:ODElin}
u'(t) \,=\, L u(t) + g(t) \,.
\eeq
Let $Z_k = \dt L_k$, corresponding to the splitting
$Z = \dt L = Z_1 + \cdots + Z_r$.
Below some rational or polynomial expressions in the matrices $Z_j$
will arise.  For this we will use the notation
\eq
\uZ \,=\, (Z_1, Z_2,\ldots, Z_r) \,.
\eeq

\subsection{Perturbed schemes}

To derive recursions for the global errors, it is
convenient to first study the effect of perturbations on the stages.
Along with (\ref{eq:PRK}) we consider a perturbed scheme
\begin{subequations}
\label{eq:PRKpert}
\setlength{\arraycolsep}{1mm}
\begin{eqnarray}
\tilde{v}_{n,i}  &  =  &  \tilde{u}_{n} + \dt \sum_{k=1}^r \sum_{j=1}^{s}
a_{ij}^{(k)} F_k(t_n + c_j \dt, \tilde{v}_{n,j})
\,+\, \rho_{n,j}  \,, \qquad
i = 1,\ldots,s\,, \qquad
\\
\tilde{u}_{n+1} &  =  &  \tilde{u}_{n} + \dt \sum_{k=1}^r \sum_{j=1}^{s}
b_{j}^{(k)} F_k(t_n + c_j \dt, \tilde{v}_{n,j})
\,+\, \sigma_{n} \,.
\end{eqnarray}
\end{subequations}
The perturbations can be used to define residual, local errors per stage.

For the vector  $c = [c_i] \in \R^s$ we denote its $j$-th power
per component as $c^j = [c_i^{\,j}]$ for $j\ge1$, with
$c^0 = e = [1,\ldots,1]^T \in \R^s$.
To make the dimensions fitting we will use the Kronecker products
$\bb{A}_k = A_k \otimes I$, $\bb{b}_k^T = b_k^T \otimes I$,
$\bb{c^j} = c^j \otimes I$ and $\bb{e} = e \otimes I$ with $m\times m$
identity matrix $I = I_{m\times m}$.
To make the notation consistent, the $ms\times ms$ identity matrix is
denoted by $\bb{I}$.
Furthermore, we let $\bb{Z}_k = I \otimes Z_k$, $Z_k = \dt L_k$,
with $I = I_{s\times s}$.

To write the difference of (\ref{eq:PRKpert}) and (\ref{eq:PRK})
in a compact form, let also $\bb{\rho}_n = [\rho_{n,i}]$
and $\bb{v}_n = [v_{n,i}]$, $\tilde{\bb{v}}_n = [\tilde v_{n,i}]
\in \R^{sm}$.  Then
\begin{subequations}
\eq
\tilde{\bb{v}}_n - \bb{v}_n \,=\, \bb{e} (\tilde u_n - u_n)
+ \sum_{k=1}^r \bb{A}_k \bb{Z}_k (\tilde{\bb{v}}_n - \bb{v}_n)
+ \bb{\rho}_n \,,
\eeq
\eq
\tilde u_{n+1} - u_{n+1} \,=\, \tilde u_n - u_n
+ \sum_{k=1}^r \bb{b}_k^T \bb{Z}_k (\tilde{\bb{v}}_n - \bb{v}_n)
+ \sigma_n \,.
\eeq
\end{subequations}
Elimination of $\tilde{\bb{v}}_n - \bb{v}_n$ leads to
\eq
\label{eq:vecform}
\tilde u_{n+1} - u_{n+1} \,=\, {R}(\uZ) (\tilde u_n - u_n)
+ \bb{r}({\uZ})^T \bb\rho_n + \sigma_n \,,
\eeq
where the amplification matrix ${R}(\uZ) \in \R^{m\times m}$ and
$\bb{r}(\uZ)^T \in \R^{m\times m s}$ are defined by
\begin{subequations}
\label{eq:Sr}
\eq
{R}(\uZ) \,=\, I + \bb{r}({\uZ})^T \bb{e} \,,
\eeq
\eq
\bb{r}({\uZ})^T \,=\, \big[{r}_1({\uZ}),\ldots,{r}_s({\uZ}) \big]
= \Big(\sum_{k=1}^r \bb{b}_k^T \bb{Z}_k \Big)
\Big(\bb{I} - \sum_{k=1}^r \bb{A}_k \bb{Z}_k \Big)^{-1}.
\eeq
\end{subequations}

The $r_j(\uZ)$ are polynomial expressions (for explicit methods) or rational
expressions (for implicit methods) in $Z_1,Z_2,\ldots,Z_r$.
It will be assumed that these expressions are bounded,
\eq
\label{eq:r_jBound}
\|r_j(\uZ)\| \, \le \, M  \qquad \mbox{for $j=1,\ldots,r$} \,,
\eeq
with some fixed $M > 0$.
For explicit methods, this will be ensured by requiring that the matrices
$Z_j$ are bounded.
Moreover, it may be assumed that the functions $r_j(\uZ)$ are not linearly 
dependent:
\eq
\label{eq:r_jIndep}
\sum_{j=1}^s \gamma_j \, r_j(\uZ) \,=\,  0 \quad 
\mbox{(for all $\uZ = (Z_1,\ldots, Z_r)$)} 
\quad \Longrightarrow \quad
\gamma_1 = \gamma_2 = \cdots = \gamma_s = 0 \,.
\eeq
Violation of this last condition would mean there are perturbations
in (\ref{eq:PRKpert}) of the form $\rho_{n,j} = \gamma_j \rho_0$, with
$\rho_0\in\R^m$, which do not influence the outcome, no matter how
large $\|\rho_0\|$ is. This indicates a redundancy (reducibility)
in the scheme.

\subsection{Error recursions}

Let $\eps_n = u(t_n) - u_n$ be the global discretization error at
time level $t_n$, $n\ge0$. As we will see, these global errors satisfy 
a recursion
\eq
\label{eq:GloErrRec}
\eps_{n+1} \,=\, {R}(\uZ) \, \eps_{n} \,+\, \delta_n \,, \qquad n \ge 0 \,,
\eeq
where $\delta_n$ is a local discretization error, introduced in the step
from $t_n$ to $t_{n+1}$.

\begin{Lem}
Suppose the functions $\fhi_k(t) = F_k(t,u(t))$ ($1\!\le\! k \!\le\! r$) 
are $l$ times continuously differentiable, and (\ref{eq:r_jBound}) is valid.
Then the local error $\delta_n$ in (\ref{eq:GloErrRec}) is given by
\eq
\label{eq:LocErr}
\delta_n \,=\, \sum_{j=1}^l \frac{\dt^j}{j!}
\sum_{k=1}^{r} {d}_{j,k}(\uZ) \fhi_k^{(j-1)}(t_n)
\,+\, \Oh(\dt^{l+1}) \max_{k, t}\|\fhi_k^{(l)}(t)\| \,,
\eeq
where
\eq
\label{eq:ErrCoef}
{d}_{j,k}(\uZ) = \big(I - j \bb{b}_k^T \bb{c}^{j-1}\big) +
\bb{r}(\uZ)^T \big(\bb{c}^j - j \bb{A}_k \bb{c}^{j-1}\big) \,.
\eeq
\end{Lem}
{\bf Proof.}
Consider the perturbed scheme (\ref{eq:PRKpert}) with $\tilde u_n = u(t_n)$
and $\tilde v_{n,i} = u(t_{n}+c_i\dt)$, $i=1,\ldots,s$.
This choice for the $\tilde v_{n,i}$ defines the perturbations $\rho_{n,i}$
and $\sigma_n$, and we obtain by Taylor expansion
$$
\begin{array}{l}
\displaystyle
\bb\rho_n = \sum_{k=1}^{r} \sum_{j\ge 1}^{} \mfrac{\dt^j}{j!}
\big(\bb{c}^j - j \bb{A}_k \bb{c}^{j-1}\big) {\fhi}_k^{(j-1)}(t_{n}) \,,
\\
\displaystyle
\sigma_n = \sum_{k=1}^{r}\sum_{j\ge 1}^{} \mfrac{\dt^j}{j!}
\big(I - j \bb{b}_k^T \bb{c}^{j-1}\big) {\fhi}_k^{(j-1)}(t_{n}) \,.
\end{array}
$$
If the $\fhi_k$ are $l$ times continuously differentiable, the sum 
over $j$ can be truncated, with $j$ ranging from $1$ to $l$ and with a 
remainder term at the $\dt^{l+1}$ level, involving the $\|\fhi^{(l)}(t)\|$
with $t$ between $t_n$ and $t_n+c_i\dt$, $i=1,\ldots,s$.

Subtraction of (\ref{eq:PRK}) from (\ref{eq:PRKpert}) 
leads to the error recursion (\ref{eq:GloErrRec}) with
$$
\delta_n = \bb{r}({\uZ})^T \bb\rho_n + \sigma_n \,.
$$
Insertion of the Taylor expansions for $\bb\rho_n$ and $\sigma_n$ thus lead
to the expressions (\ref{eq:LocErr}), (\ref{eq:ErrCoef}) for the local 
errors.
\hfill $\Box$

\bigskip
For a method with classical order $p$ and stage order $q\le p$, we have
\begin{subequations}
\eq
{d}_{j,k}(\uZ) \,=\, 0  \qquad\mbox{if $\; j \le q$} \,,
\eeq
\eq
{d}_{j,k}(\uZ) \,=\, 
\bb{r}(\uZ)^T \big(\bb{c}^j - j \bb{A}_k \bb{c}^{j-1}\big)
\qquad \mbox{if $\; q < j \le p$} \,.
\eeq
\end{subequations}
Note that $\bb{r}(\uZ) = 0$ if all $Z_l=0$, and so the same property
holds for the functions ${d}_{j,k}(\uZ)$, $q < j \le p$.  In fact,
since we know that $\delta_n = \Oh(\dt^{p+1})$ for non-stiff problems,
it follows that $d_{j,k}(\uZ) = \Oh(\dt^{p+1-j})$ if all $Z_l = \Oh(\dt)$,
$1\le l\le r$.
The above properties will be used in the analysis of the local
discretization errors.

\begin{Rem} \rm \label{Rem:Nonlin}
The above derivations can also be performed for nonlinear problems
(\ref{eq:ODE}), essentially by replacing occurring Kronecker products 
such as $\bb{Z} = I\otimes Z$, $Z = \dt L$,
with the varying block-diagonal matrix
$\bb{Z} = \mbox{Diag}(Z_{n,i}) \in \R^{ms\times ms}$
where
${Z}_{n,i} (\tilde v_{n,i} - v_{n,i}) =
\dt \big(F(\tilde v_{n,i}) - F(v_{n,i})\big)$,
with changes over the steps and the stages.
However, this leads to more complicated notation, and it does
not give additional insight.
\hfill $\Diamond$
\end{Rem}

\section{Error analysis for cell-based splittings}
\label{sect:cell-based}

From now on, we restrict our attention to explicit methods.
In this section it will be assumed that the splitting (\ref{eq:split})
is cell-based, $F_j = I_j F$ for $j=1,\ldots,r$. Then we have
$\fhi_k(t) = I_k u'(t)$, which is bounded in the maximum norm
uniformly in the spatial mesh width. For flux-based splittings,
considered in Section~\ref{Sect:flux-based},
this last property will not be valid.

Throughout the remaining sections we will denote by $\Oh(\dt^q)$ a scalar 
or vector for which all components can be bounded $K \dt^q$, for $\dt > 0$ 
small enough, with $K$ not depending on the mesh widths $\dx_j$ in the spatial
discretization.

\subsection{Order reduction}

In this section we derive bounds for the discretization errors that
are valid for semi-discrete systems with smooth solutions.
The classical, non-stiff order conditions are then no longer sufficient
to obtain convergence of order $p$. This so-called order reduction is
due to the fact that $F$ contains
negative powers of the mesh widths $\dx_j$ in space.
We will accept a restriction on $\dt/\dx_j$ for stability, but the
resulting error bounds should not contain negative powers of $\dx_j$.

For the partitioned methods we want to see the effects of the partitioning
on the errors. We will therefore study the errors in the maximum norm,
assuming stability of the scheme:
\eq
\label{eq:stab}
\sup_{n\ge0}\|{R}(\uZ)^n\|_\infty \le K \,.
\eeq
Sufficient conditions for having (\ref{eq:stab}) with $K=1$ have been derived 
in \cite{HMS13} for nonlinear problems.
For explicit methods, a necessary stability condition is boundedness of the 
$Z_j$, and therefore (\ref{eq:r_jBound}) will be satisfied.

Let $Z_k = \dt L_k$, corresponding to the splitting
$L = L_1 + \cdots + L_r$.
If $L$ is a discretized convection operator, and a CFL restriction
$\dt/\dx_j \le \nu$ is satisfied with some fixed $\nu$, then
$\|Z_k\|_\infty = \Oh(1)$.
It will be tacitly assumed that the exact solution is smooth,
so that derivatives of $u(t)$ are $\Oh(1)$.
If the splitting is cell-based, then $\fhi_k(t) = I_k u'(t)$, so any term
$\fhi_k(t)$ and its time derivatives will then be $\Oh(1)$. Note, however,
that $\fhi_k(t)$ is not a smooth grid function: there will be jumps over
the interfaces of the spatial components, and therefore we will in general
only have $\|Z \fhi_k(t)\|_\infty = \Oh(1)$ instead of 
$\|Z \fhi_k(t)\|_\infty = \Oh(\dt)$.

If the stability assumption (\ref{eq:stab}) holds, it follows directly
that consistency of order $q$ (i.e., $\|\delta_n\|_\infty = \Oh(\dt^{q+1})$)
implies convergence of order $q$ (i.e., $\|\eps_n\|_\infty = \Oh(\dt^q)$),
but we will see that the order of convergence can also be one larger than
the order of consistency.

\subsection{Local error analysis}

To analyze the order of the local errors we will distinguish various
cases , depending whether the method is internally consistent or not
(stage order $q\ge1$ or $q=0$).

{\it Stage order zero\/}:
Let us first consider a method with classical order $p\ge 1$ and stage order
$q = 0$, that is, the method is not internally consistent:
${A}_k e \neq {A}_l e$  for some $k,l$.
Then the leading term in the local error is
\eq
\label{eq:locerr1}
\delta_n = \dt\;  \sum_{k=1}^r {d}_{1,k}(\uZ)\fhi_k(t_n) + \Oh(\dt^2) \,,
\qquad
{d}_{1,k}(\uZ) = \bb{r}(\uZ)^T (\bb{c} - \bb{A}_k \bb{e}) \,.
\eeq
Since $\fhi_k(t_n) = \Oh(1)$, this gives an $\Oh(\dt)$ local error bound
in the maximum norm, which is of course quite poor. 
After all, $\delta_n$ is the error that results after one step if $\eps_n = 0$.
However, it will be seen that this still can lead to convergence of order one.

{\it Stage order one\/}:
Next assume $q\ge1$, that is, the internal consistency condition 
(\ref{eq:intcons}) is satisfied:  $A_k e = A_l e$ for $1\le k,l\le r$.
If $p = 1$ it follows directly that $\|\delta_n\|_\infty = \Oh(\dt^2)$.
If $p\ge2$ the leading term in the local discretization errors is
given by
\eq
\label{eq:locerr2}
\delta_n = \sfrac{1}{2} \dt^2 \, \sum_{k=1}^r
{d}_{2,k}(\uZ) \fhi_k'(t_{n}) + \Oh(\dt^3) \,,
\qquad
{d}_{2,k}(\uZ) = \bb{r}(\uZ)^T
\big( \bb{c}^2 - 2 \bb{A}_k \bb{c}\big)
\,.
\eeq
This still gives only consistency of order one, that is, an error
$\Oh(\dt^2)$ after one step,
but we will discuss below damping and cancellation effects that can
lead to convergence with order two in this case.

{\it Higher orders\/}:
For explicit methods it is not possible to have $c^2 = 2 A_k c$. With 
(\ref{eq:r_jIndep}) this implies that the functions ${d}_{2,k}(\uZ)$ 
cannot be identically equal to zero.
Yet there are exceptional cases where (\ref{eq:locerr2}) can give
consistency of order larger than one, under the assumption that
$L u''(t) = \Oh(1)$.
If all ${d}_{2,k}(\uZ)$ are equal, say
\eq
\label{eq:dsequal1}
{d}_{2,k}(\uZ) \,=\, Q(\uZ) \qquad \mbox{for $\;k=1,\ldots,r$} \,,
\eeq
with $Q(\uZ)$ a polynomial expression in the $Z_k$,
then we have $\delta_n = \frac{1}{2} \dt^2 Q(\uZ) u''(t_n) + \Oh(\dt^3)$.
Since the constant term in $Q(\uZ)$ is zero if $p\ge2$, and 
$Z_k u''(t) = \dt I_k L u''(t) = \Oh(\dt)$,
this gives indeed $\|\delta_n\|_\infty = \Oh(\dt^3)$.
It should be noted, however, that (\ref{eq:dsequal1}) will only occur if
\eq
\label{eq:dsequal2}
A_k c = A_l c \qquad \mbox{for all $\;k,l=1,\ldots,r$} \,.
\eeq
This will hold, of course, for the case that all coefficient matrices $A_k$ 
are equal, but it is not valid for the multirate methods from 
Section~\ref{sect:Examples}.

Methods with equal coefficient matrices $A_k$ have been studied in
\cite{KMR13}. For such methods the above arguments can be simplified,
see Theorem~2.1 in \cite{KMR13} and Remark~\ref{Rem:EqualA}
in the present paper, since the internal stages then only 
use the complete function $F$ rather than the $F_k$ from the decomposition
(\ref{eq:split}).

The above expressions for the local errors show that {\it order reduction} 
is to be expected: the accuracy will primarily depend on the stage order
$q$, rather than on the classical order $p$.
This order reduction will appear primarily at interface points on the
spatial grid, where the grid-functions $\fhi_k(t)$ have jumps.

Further we note that these expressions for the local errors are similar
to those given, for example, in \cite{HuRu07} for implicit-explicit 
Runge-Kutta methods, and in \cite{PBJ04a} for a class of implicit 
additive Runge-Kutta methods for parabolic problems with domain decomposition.

\subsection{Global error analysis}

Based on the local error behaviour, one would expect convergence
with order one for the TW2 and SH2 schemes, and lack of convergence
for the scheme CS2. This is not what was seen in the numerical test
for advection with a smooth solution.
To obtain the correct (observed) order of convergence, we
need to study the propagation of the leading term in the local error.

In the following result we consider a partitioned method (\ref{eq:PRK}) 
with classical order $p$ and stage order $q$. For the leading local 
error terms, it will be assumed that there is a matrix $W \in \R^{m\times m}$ 
such that
\eq
\label{eq:Wbound}
\big(\bb{r}(\uZ)^T \bb{e}\big) W  \,=\, \sum_{k=1}^r {d}_{q+1,k}(\uZ) I_k
\,.
\eeq

\begin{Thm}
\label{Thm:conver}
Assume that $p \ge q+1$ and the stability condition (\ref{eq:stab})
holds. Assume furthermore that (\ref{eq:Wbound}) holds with
$\|W\|_\infty = \Oh(1)$.
Then the method is convergent with order $q+1$ in the maximum norm.
\end{Thm}

\noindent
{\bf Proof.}
Let $\xi_n = \frac{1}{(q+1)!} \dt^{q+1} W u^{(q+1)}(t_n)$. Then
$\|\xi_n\|_\infty = \Oh(\dt^{q+1})$,
$\|\xi_{n+1} - \xi_n\|_\infty = \Oh(\dt^{q+2})$,
and the local error $\delta_n$ from (\ref{eq:LocErr}) can be
decomposed as
$$
\delta_n \,=\, (R(\uZ) - I) \xi_n + \eta_n \,,
$$
where $\eta_n = \Oh(\dt^{q+2})$ contains the higher-order terms.
Introducing $\hat{\eps}_n = \eps_n - \xi_n$, we get the recursion
$$
\hat{\eps}_{n+1} \,=\, R(\uZ) \hat{\eps}_n + \hat{\delta}_n \,, \qquad
\hat{\delta}_n \,=\, \xi_{n+1} - \xi_n + \eta_n \,.
$$
In the standard way, it is seen from (\ref{eq:stab}) that
$
\| \hat{\eps}_n\|_\infty \le K ( \| \hat{\eps}_0\|_\infty +
\sum_{j=0}^n \| \hat{\delta}_j \|_\infty ).
$
Since $\eps_0 = 0$ we obtain
$$
\|\eps_n\|_\infty \le \|\xi_n\|_\infty + K\|\xi_0\|_\infty +
\sum_{k=0}^n K \big( \|\xi_{k+1} - \xi_k\|_\infty + \|\eta_k\|_\infty \big) \,,
$$
from which the convergence result now follows.
\hfill $\Box$

\bigskip
This result and its proof is similar as for standard Runge-Kutta methods
where order reduction may arise due to boundary conditions; see e.g.\
\cite{BCT82} or the review in \cite[Sect.\,II.2]{HuVe03}.
With the partitioned Runge-Kutta methods and  multirate schemes, we are 
creating interfaces that act like (internal) boundaries with time-dependent
boundary conditions.

\subsection{Examples: multirate methods with cell-based splittings}
\label{subsect:MR2}

For the simple multirate examples from Section~\ref{sect:MR1} we
will study the order of convergence in the maximum norm. It will be
assumed that
\eq
\label{eq:stabcond}
\|I + Z_1\|_\infty \,\le 1\, \,, \qquad
\|I + \sfrac{1}{2}Z_2\|_\infty \,\le\, 1 \,.
\eeq
These conditions (or, rather, the nonlinear counterparts) were used
in \cite{HMS13} to prove the stability condition (\ref{eq:stab})
with $K=1$ for the multirate schemes.

\subsubsection{First-order multirate schemes OS1, TW1}

For the TW1 scheme, we have $p = q = 1$, giving local errors
$\|\delta_n\|_\infty = \Oh(\dt^2)$ from which we obtain in the
standard way convergence with order $1$ in the maximum norm.

Consider the OS1 scheme, with $p=1$ but $q=0$. Here
$\|\delta_n\|_\infty = \Oh(\dt)$ only. Still, first-order convergence
can be shown. For this, it is assumed, in addition to (\ref{eq:stabcond}),
that
\eq
\label{eq:stabcond'}
\|Z_2\|_\infty \,\le\, 4\theta \,<\, 4 \,.
\eeq
From (\ref{eq:stabcond}) it follows already that $\|Z_2\|_\infty \le 4$,
and consequently (\ref{eq:stabcond'}) is only a minor strengthening of
the assumptions (\ref{eq:stabcond}).

For this OS1 scheme we have, with $Z = Z_1 + Z_2$,
$$
\bb{r}(\uZ)^T = \big[\sfrac{1}{2} Z (I + \sfrac{1}{2} Z_2) \quad
\sfrac{1}{2} Z\big] \,,
\qquad
c - A_1 e = \left[ \begin{array}{c}
0 \\ \sfrac{1}{2} 
\end{array} \right] \,,
\qquad
c - A_2 e = \left[ \begin{array}{c}
0 \\ 0
\end{array} \right]  \,,
$$
Hence
$$
\bb{r}(\uZ)^T \bb{e} = Z (I + \sfrac{1}{4} Z_2) \,, \qquad
{d}_{1,1}(\uZ) = \sfrac{1}{4} Z \,, \qquad
{d}_{1,2}(\uZ) = 0 \,,
$$
and (\ref{eq:Wbound}) reads
$$
(I + \sfrac{1}{4} Z_2) W = I_1 \,.
$$
In view of (\ref{eq:stabcond'}) we have
$\|(I + \sfrac{1}{4} Z_2)^{-1}\|_\infty \le (1 - \theta)^{-1}$,
which ensures the bound $\|W\|_\infty \le (1 - \theta)^{-1}$.
Application of Theorem~\ref{Thm:conver} shows that the OS1 scheme
will indeed converge with order $1$ under the assumptions
(\ref{eq:stabcond}), (\ref{eq:stabcond'}).

\subsubsection{Second-order multirate schemes CS2, TW2, SH2}

For the second-order methods, the expressions for $R(\uZ)$ and the 
$d_{q+1,k}(\uZ)$ are already rather complicated. Therefore, instead of 
a detailed analysis of (\ref{eq:Wbound}), we will only present 
here some experimental results for the semi-discrete system
\eq
\label{eq:adveq_sd}
u'_j(t) = \mfrac{1}{\dx_j} \big(u_{j-1}(t) - u_j(t)\big)
\qquad\mbox{for $j\in\I = \{1,2,\ldots,m\}$} \,,
\eeq
with $u_0(t)=0$, corresponding to first-order upwind discretization of
the advection equation $u_t + u_x = 0$ with homogeneous inflow condition
$u(0,t) = 0$.

\begin{figure}[htb]
\setlength{\unitlength}{1cm}
\begin{center}
\begin{picture}(13,3.5)
\epsfxsize13cm
\put(0,0){\epsfbox{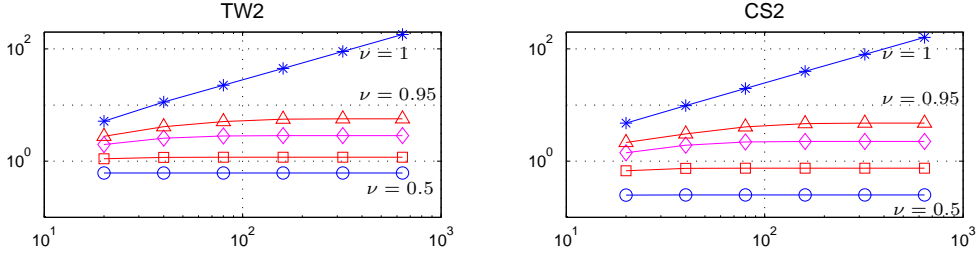}}
\put(4.80, 0.7){\scriptsize $\nu = 0.5$}
\put(11.8,0.43){\scriptsize $\nu = 0.5$}
\put(4.70, 1.95){\scriptsize $\nu = 0.95$}
\put(11.65,1.9){\scriptsize $\nu = 0.95$}
\put(4.7, 2.5){\scriptsize $\nu = 1$}
\put(11.65,2.5){\scriptsize $\nu = 1$}
\end{picture}
\caption{ \small \label{Fig:boundW}
Norm $\|W\|_\infty$ versus $m = 20, 40,\ldots,640$ for various
values of $\nu = \dt/h$ with the schemes TW2 (left) and CS2 (right).
Markers: {\large $\circ$} for $\nu = 0.5$,
{\tiny $\square$} for $\nu = 0.75$, $\diamond$ for $\nu = 0.9$,
{\tiny $\bigtriangleup$} for $\nu = 0.95$ and $*$ for $\nu = 1$.
}
\vspace{-2mm}
\end{center}
\end{figure}

We take a partitioning $\I = \I_1 \cup \I_2 = \{1,2,\ldots,m\}$ with
$\I_2 = \{j : {1\over4}m < j \le {3\over4}m\}$,
and mesh widths $\dx_j = h$ if $j\in\I_1$,
$\dx_j = {1\over2} h$ if $j\in\I_2$, with $h = 4/(3 m)$.
In Figure~\ref{Fig:boundW} the norm $\|W\|_\infty$ is plotted as
function of $m=20, 40, \ldots, 640$ for various values of $\nu = \dt/h$
for the schemes TW2 and CS2; the results for SH2 were similar to
those of TW2. In this example, the matrix ${r}(\uZ)^T \bb{e}$ is
nonsingular, and it is well-conditioned for $\nu \le 1$.
We see that $\|W\|_\infty = \Oh(1)$ provided that $\nu < 1$,
whereas $\|W\|_\infty \sim m$ if $\nu=1$.
Other partitionings $\I = \I_1 \cup \I_2$ produced similar results.

The combination of Theorem~\ref{Thm:conver} and these experimental
bounds for first-order advection discretization does provide a heuristic
explanation for the advection test results in Section~\ref{Ssect:numtests1},
where we observed convergence of the schemes TW2 and SH2 with order two in
the maximum norm, and with order one for the CS2 scheme.

\subsection{Numerical test: 2D advection}
\label{Ssect:numtests}

The numerical test in Section~\ref{Ssect:numtests1} for 1D advection was
highly artificial, because there was no practical need to refine on
subintervals. Below we will present a more relevant test for advection in 2D.

As before, we will use the WENO5 scheme for the spatial discretization.
This spatial scheme combines high accuracy with a good behaviour near
discontinuities. Since the focus in this paper is temporal accuracy,
we will use linear advection problems with smooth initial profiles
in the tests. Due to the WENO5 spatial discretization, the semi-discrete
ODE system is still nonlinear.

As a test example we consider here the two-dimensional advection equation
\begin{subequations}
\label{eq:Adv2D}
\eq
u_t + (a_1 u)_x + (a_2 u)_y \,=\, 0
\eeq
for $0<x,y<1$, $0<t\le1$, with divergence-free velocity field given by
\eq
a_1(x,y) = 2 \pi (y-\sfrac{1}{2}) \,, \qquad 
a_2(x,y) = -2 \pi (x-\sfrac{1}{2}) \,, 
\eeq
and initial profile 
\eq
u(x,y,0) \,=\, e^{-10\big((x-\frac{1}{2})^2+(y-\frac{1}{4})^2\big)} \,.
\eeq
\end{subequations}
The wind field gives a uniform clock-wise rotation around the center of 
the domain.  We take end time $T=\frac{1}{3}$, giving a rotation of the 
initial profile over an angle $\frac{2}{3}\pi$.  
At the inflow boundaries Dirichlet conditions are described, corresponding 
to the exact solution.

\begin{figure}[htb]
\setlength{\unitlength}{1cm}
\begin{center}
\begin{picture}(13,5.0)
\epsfxsize13cm
\put(0,0){\epsfbox{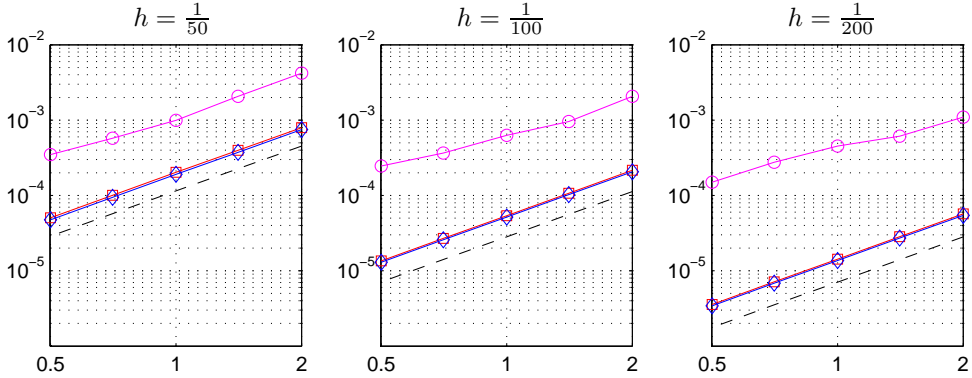}}
\put(1.7 ,4.6){$h = \frac{1}{50}$}
\put(6.0 ,4.6){$h = \frac{1}{100}$}
\put(10.4,4.6){$h = \frac{1}{200}$}
\end{picture}
\caption{ \small \label{Fig:A2Dres00}
Test (\ref{eq:Adv2D}) with cell-based decomposition.
Maximum errors for the schemes
TW2 ({\tiny $\square$} marks), CS2 ($\circ$ marks) and SH2 ($\diamond$ marks)
as function of Courant numbers $\nu =  2\pi\dt/h$, between $\frac{1}{2}$ 
and $2$,
for the grids with $h = \frac{1}{50}$ (left), $h = \frac{1}{100}$ (middle)
and $h = \frac{1}{200}$ (right).
}
\vspace{-2mm}
\end{center}
\end{figure}

We consider a partitioning where ${\cal I}_1$ corresponds to the grid points 
in the region where $|x-\frac{1}{2}| + |y-\frac{1}{2}| \le \frac{1}{3}$.
This is a natural partitioning since the velocity field increases
towards the corners of the domain.
In the test we compare the solutions obtained by the multirate
schemes with an accurate semi-discrete solution, obtained with a
Runge-Kutta method with small step-size.

The results on three uniform grids, with  $\dx = \dy = h$, $h = \frac{1}{50}$,
$\frac{1}{100}$, $\frac{1}{200}$,
are presented in Figure~\ref{Fig:A2Dres00}. There, for each separate grid,
the maximum errors are plotted for various Courant numbers
$\nu = \dt\max_{x,y}(|a_1|+|a_2|)/h = 2\pi\dt/h$.
The dashed line in the figures gives the result for the scheme
where in each time step the explicit trapezoidal rule is applied twice,
with step-size $\frac{1}{2}\dt$, over the whole region.

On any fixed grid the three schemes are second-order convergent (classical
order two), but it is clear that the CS2 scheme has a large error constant,
affected by $h$.
Comparing the errors on the three grids for the same Courant number shows
indeed a very slow convergence for the CS2 scheme.

\section{Decomposition based on fluxes}
\label{Sect:flux-based}

For conservation laws, the semi-discrete system (\ref{eq:ODE}) will in
general be in conservative form. In 1D, for example, we will have
\eq
\label{eq:Cons1D}
u'_j(t) =
\mfrac{1}{\dx_j} \big(f_{j-\half}(u(t)) - f_{j+\half}(u(t)) \big)
\,, \qquad j \in \I = \{1,2,\ldots,m\} \,.
\eeq
Multirate methods can be based on these numerical fluxes $f_{j\pm1/2}(u)$
rather than in terms of the components of $F(u)$.

A decomposition $F = F_1 + F_2 + \cdots + F_r$ can be based on
fluxes in the following way.
The conservative semi-discrete ODE system (\ref{eq:Cons1D})
has right-hand side function
\eq
F(v) = H^{-1} D \, \Phi(v)
\eeq
with $H = \mbox{diag}(\dx_j)$, bi-diagonal difference matrix $D$, 
and flux vector $\Phi(v) = [\Phi_j(v)]$, $\Phi_j(v) = f_{j+1/2}(v)$.
If $J_k$ corresponds to a discrete indicator function for a region
$\Omega_k$ of the PDE domain $\Omega = \Omega_1\cup\ldots\cup\Omega_r$, 
then 
\eq
F_k(v) = H^{-1} D \, J_k \, \Phi(v) \,, \qquad k=1,\ldots,r\,,
\eeq
gives a flux-based decomposition of $F$.

As an example, suppose that $r=2$, $\Omega_1 = \{x : x \le x_i\}$ and 
$\Omega_2 = \{x : x \ge x_i\}$. Then the $j$th component of 
the vector functions $F_{1}$ and $F_{2}$ is given by
\eq
\label{eq:Exa1Dsplit}
\left. \begin{array}{ccl}
F_{1,j}(v) = \mfrac{1}{\dx_j} \big(f_{j-\half}(v) - f_{j+\half}(v) \big) \,,
\qquad \qquad
F_{2,j}(v) = 0
& &  \mbox{for $j<i$} \,,
\\[2mm]
F_{1,j}(v) = \mfrac{1}{\dx_i}f_{i-\half}(v) \,,
\qquad \qquad \quad
F_{2,j}(v) = \mfrac{-1}{\dx_i}f_{i+\half}(v)
& &  \mbox{for $j=i$} \,,
\\[2mm]
F_{1,j}(v) = 0 \,,
\qquad \qquad
F_{2,j}(v) = \mfrac{1}{\dx_j} \big(f_{j-\half}(v) - f_{j+\half}(v) \big)
& & \mbox{for $j>i$} \,.
\end{array}
\right\}
\eeq
Since we are here dealing with fluxes, mass-conservation is guaranteed
for any stage.  However, there some serious drawbacks as well.

First, monotonicity assumptions such as $\|v+\tau F_{k}(v)\|\le\|v\|$
will not be valid in the maximum norm with this decomposition.
This can be seen already quite easily for the first-order upwind
advection discretization (\ref{eq:adveq_sd}) with $r=2$.
Writing this system as $u'(t) = L u(t)$, the above decomposition
would correspond to $L = L I_1 + L I_2$, that is, $F_{k} = L I_k$,
but it is easy to show that $\|I + \tau L I_k\|_\infty$ is larger than
one for any $\tau > 0$.
Consequently, stability assumptions like (\ref{eq:stabcond}) are also no 
longer relevant.

Secondly, such a flux-based decomposition of $F$ can easily lead to 
inconsistencies, since we do not have $F_{k}(u(t)) = \Oh(1)$, no matter 
how smooth the solution is.
For example, for the first-order upwind system (\ref{eq:adveq_sd}),
using these $F_{1}$ and $F_{2}$ in the OS1 scheme gives a completely
inconsistent scheme.
This issue of accuracy will be discussed next.

\subsection{Error analysis}

We will discuss here the effect of flux-based decompositions on the local 
errors. The transition of local to global errors is similar to the
cell-based decompositions.
Note that the formulas (\ref{eq:locerr1}) and (\ref{eq:locerr2}) are 
still correct for the leading term, with $\fhi_k(t) = F_k(t,u(t))$. 
However, now $\|\fhi_k(t)\|_\infty$ will be proportional to $1/\dx$,
see e.g.\ formula (\ref{eq:Exa1Dsplit}) with $j=i$, and therefore we only 
have $\|\dt\fhi_k(t)\|_\infty = \Oh(1)$ under a CFL restriction on $\dt/\dx$.
This may lead to smaller orders of consistency/convergence than for the 
cell-based splittings.
We will discuss various cases, leading to convergence with order zero, 
one, or more, separately.

{\it Stage order zero\/}:
If the method is not internally consistent, the principal local error term
is $\dt \sum_k d_{1,k}(\uZ) \fhi_k(t_n)$, see (\ref{eq:locerr1}). Since
we now only have $\dt\fhi_k(t) = \Oh(1)$, the error after one step may
not tend to zero as $\dt\rightarrow0$.

\begin{Exa} \rm
For the advection equation $u_t + u_x = 0$, consider (\ref{eq:Exa1Dsplit}) 
with first-order upwind fluxes $f_{j+1/2}(v) = v_j$,
and denote the components of the vector $u_n$ as 
$u_j^n \approx u(x_j,t_n)$. A little calculation shows that at the 
interface point the scheme (\ref{eq:OS1}) gives
$$
u_i^{n+1} = u_i^n + \mfrac{\dt}{\dx_i} (u_{i-1}^n - u_i^n)
+ \sfrac{1}{4} \Big(\mfrac{\dt}{\dx_i}\Big)^2 u_i^n \,.
$$
Already after one step, starting with $u_j^0 = u(x_j,t_0)$, 
this gives an error 
$u_i^1 - u(x_i,t_1) = \frac{1}{4} \nu_i^2 u(x_i,t_0) + \Oh(\dt)$ if
$\dt\rightarrow0$ while $\nu_i = \dt/\dx_i$ is held fixed, leading to
an $\Oh(1)$ error in the maximum norm if $u(x_i,t_0) \neq 0$.
\hfill $\Diamond$
\end{Exa}

{\it Stage order one\/}: 
For methods that are internally consistent, with stage order $q=1$, the
the principal local error term is 
$\frac{1}{2}\dt^2 \sum_k d_{2,k}(\uZ) \fhi'_k(t_n)$, see (\ref{eq:locerr2}). 
Since $\dt\fhi'_k(t) = \Oh(1)$, this gives an error proportional to $\dt$
in the maximum norm after one step.
Due to damping and cancellation effects, we can still have convergence 
with order~$1$. In some numerical tests this will be seen to hold for 
the multirate methods from Section~\ref{sect:Examples}.

{\it Higher orders\/}:
If we have an internally consistent method, for which all $d_{2,k}(\uZ)$ 
are equal, say ${d}_{2,k}(\uZ) = Q(\uZ)$ for $k=1,\ldots,r$,
as in (\ref{eq:dsequal1}), then $\|\delta_n\|_\infty = \Oh(\dt^2)$, 
because $\sum_k \fhi_k(t) = u'(t)$, which is a smooth, bounded grid 
function, unlike the individual $\fhi_k(t)$ terms.
As noted before, this requires (\ref{eq:dsequal2}), which does not
hold for the multirate methods from Section~\ref{sect:Examples}.

For general partitioned methods, if we have, instead of (\ref{eq:dsequal1}), 
the stronger assumption
\eq
\label{eq:dsequal2'}
{d}_{2,k}(\uZ) \,=\, P(\uZ)\cdot Z \qquad \mbox{for $k=1,\ldots,r$} \,,
\eeq
then the leading term in $\|\delta_n\|_\infty$ will even be $\Oh(\dt^3)$, 
because in this case
$$
\sfrac{1}{2}\dt^2 \sum_{k=1}^r d_{2,k}(\uZ) \fhi'_k(t_n) =
\sfrac{1}{2}\dt^2 P(\uZ)\, Z u''(t_n) \,,
$$ 
and $Z u''(t) = \dt L u''(t) = \Oh(\dt)$ if $L u''(t) = \Oh(1)$, which
will be valid if the PDE solution is smooth with boundary conditions that
are constant in time.
The assumption (\ref{eq:dsequal2'}) will hold if $p\ge3$ and all 
coefficient matrices $A_k$ are equal.

\begin{Rem} \rm \label{Rem:EqualA}
As noted above, having a partitioned method with equal coefficient
matrices $A_k$ will often be beneficial with respect to the accuracy.
In fact, it was shown in \cite{KMR13} that for such methods the order 
of consistency will be $p$ for cell-based splittings and $p-1$ for 
flux-based splittings.
This can also be demonstrated from the local error expansions that are
used in this paper.  

If $A_k = A$ for all $k$, then leading term in the local error is
given by
$$
\delta_n = \sfrac{1}{2} \dt^2 Q(\uZ) u''(t_n) + \Oh(\dt^3) \,,
\qquad
Q(\uZ)  \,=\, \bb{r}({\uZ})^T (\bb{c}^2 - 2 \bb{A} \bb{c}) \,,
$$
and we have
$$
\textstyle
\bb{r}({\uZ})^T \,=\, (\sum_{k} \bb{b}_k^T \bb{Z}_k )
(\bb{I} - \bb{A} \bb{Z} )^{-1}
\\
=\, (\sum_{k} \bb{b}_k^T \bb{Z}_k )
(\bb{I} + \bb{A} \bb{Z} +  \bb{A}^2 \bb{Z}^2 +\cdots ) \,.
$$
Hence
$$
\textstyle
Q(\uZ)  \,=\, \sum_{k} Z_k (q_{1k} + q_{2k} Z + q_{3k} Z^2
+ \cdots ) \,, \qquad
q_{jk} = b_k^T A^{j-1} (c^2 - 2 A c) \,.
$$
If the method has order $p$ we have $q_{jk} = 0$ for $0\le j\le p-2$,
$1\le k\le r$. Therefore, neglecting the higher order terms,
\eq
\label{eq:locerrA}
\textstyle
\delta_n \,=\, \sfrac{1}{2} \dt ^p \sum_{k} Z_k
(q_{p-1,k} + q_{p,k} Z + \cdots ) L^{p-2} u''(t_n) \,.
\eeq
Assuming $L^{p-2} u''(t) = \Oh(1)$, which is an assumption on the
boundary conditions for the PDE solution, this gives
$\delta_n = \Oh(\dt^p)$. For cell-based splitting, $Z_k = I_k L$, we get
\eq
\label{eq:locerrB}
\textstyle
\delta_n \,=\, \sfrac{1}{2} \dt ^{p+1} \sum_{k} I_k
(q_{p-1,k} + q_{p,k} Z + \cdots ) L^{p-1}  u''(t_n) \,.
\eeq
Under the assumption $L^{p-1} u''(t) = \Oh(1)$ it now follows that
$\delta_n = \Oh(\dt^{p+1})$, which is the classical order
of consistency.
\hfill $\Diamond$
\end{Rem}

\subsection{Numerical tests}
\label{Ssect:numtests2}

To show the effect of flux-based decompositions, the previous tests
are repeated with the CS2, TW2 and SH2 schemes.
It should be noted that, since the CS2 scheme is always conservative, 
there is actually no need to apply this scheme with such flux-based 
decompositions.  Instead, the more accurate cell-based decompositions 
can be used for this scheme.

\subsubsection{1D advection}

We consider once more the simple problem $u_t + u_x = 0$ with periodic 
boundary conditions and $u(x,0) = \sin^2(\pi x)$, that was already used
in Section~\ref{Ssect:numtests1} with cell-based splittings. 
The set-up of the test is the same as before, with WENO5 spatial 
discretization, temporal refinement on the domain
$\{ x : |x-\frac{1}{4}| \le \frac{1}{8}\} \cup
\{ x : |x-\frac{3}{4}| \le \frac{1}{8}\}$,
and a fixed Courant number $\nu = \dt/\dx = 0.5$,
only now we consider a flux-based splitting of $F = F_1 + F_2$.
The results are given in Table~\ref{Tab:Test1FluxBased}.

\begin{table}[htb]
\caption{ \small \label{Tab:Test1FluxBased}
Flux-based splitting.
Results for the smooth advection problem with the CS2, TW2 and SH2 schemes.
Maximum errors and $L_1$-errors at final time $t_N=1$
for various $m$ with fixed Courant number $\dt/\dx = 0.5$, $\dx = 1/m$.
The approximate order of convergence is also given.
}
\small
\begin{center}
\setlength{\tabcolsep}{3mm}
\begin{tabular}{|l||c|c|c|c||c|}
\hline \rule[-1mm]{0mm}{5mm}
\hfill $m$  $\qquad$  & 100    & 200     & 400     & 800 & Order
\\ \hline \hline
\rule{0cm}{4mm}
CS2, $\|\eps_N\|_\infty$
 & $3.98\cdot10^{-2}$ & $3.65\cdot10^{-2}$
 & $3.54\cdot10^{-2}$ & $3.52\cdot10^{-2}$
& 0 \\ \rule{0cm}{3mm}
CS2, $\|\eps_N\|_1$
 & $4.43\cdot10^{-3}$ & $1.48\cdot10^{-3}$
 & $5.12\cdot10^{-4}$ & $2.09\cdot10^{-4}$
& 1 \\[1mm] \hline
\rule{0cm}{4mm}
TW2, $\|\eps_N\|_\infty$
 & $8.20\cdot10^{-4}$ & $4.20\cdot10^{-4}$
 & $2.45\cdot10^{-4}$ & $1.31\cdot10^{-4}$
\rule{0cm}{4mm}
& 1 \\ \rule{0cm}{3mm}
TW2, $\|\eps_N\|_1$
 & $2.45\cdot10^{-4}$ & $6.57\cdot10^{-5}$
 & $1.80\cdot10^{-5}$ & $5.08\cdot10^{-6}$
& 2 \\[1mm] \hline
\rule{0cm}{4mm}
SH2, $\|\eps_N\|_\infty$
 & $3.73\cdot10^{-4}$ & $1.30\cdot10^{-4}$
 & $6.69\cdot10^{-5}$ & $3.77\cdot10^{-5}$
& 1 \\ \rule{0cm}{3mm}
SH2, $\|\eps_N\|_1$
 & $2.07\cdot10^{-4}$ & $5.29\cdot10^{-5}$
 & $1.36\cdot10^{-5}$ & $3.49\cdot10^{-6}$
& 2 \\[1mm] \hline
\end{tabular}
\end{center}
\end{table}

From this table, we make the following observations. In the maximum norm 
there is no convergence for the CS2 scheme, and only first-order 
(approximately) convergence for the TW2 and SH2 schemes.
In the $L_1$-norm these orders of convergence are one higher, due to
the fact that the largest errors are confined to relatively small
spatial regions, near the interface points.

\subsubsection{2D advection}

Also the test for 2D advection (\ref{eq:Adv2D}) with rotational velocity 
field was performed again, but now with flux-based decomposition.  
Similar as before, for the cell-based splittings, we consider a 
partitioning  of the region with step-size $\dt$ in that part of the 
region where $|x-\frac{1}{2}| + |y-\frac{1}{2}| \le \frac{1}{3}$, and 
a step-size $\frac{1}{2}\dt$ is taken elsewhere.
The errors in the solutions obtained by the multirate schemes are measured
with respect to an accurate numerical solution of the semi-discrete system,
so it is only the temporal error that is measured here. 

\begin{figure}[htb]
\setlength{\unitlength}{1cm}
\begin{center}
\begin{picture}(13,5.0)
\epsfxsize13cm
\put(0,0){\epsfbox{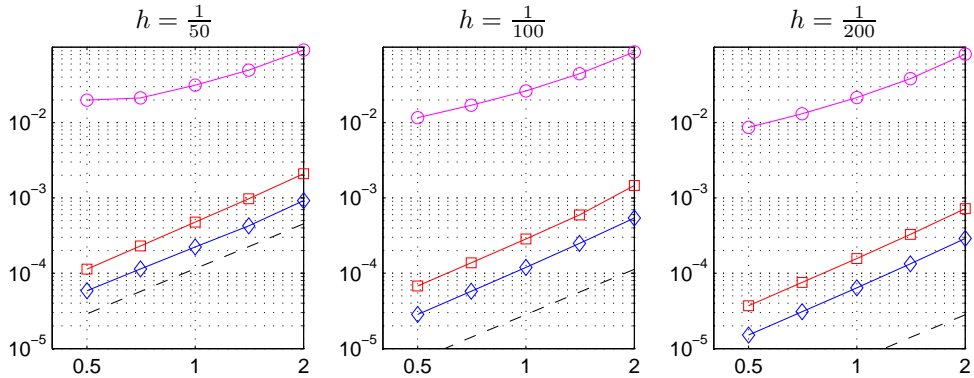}}
\put(1.7 ,4.6){$h = \frac{1}{50}$}
\put(6.0 ,4.6){$h = \frac{1}{100}$}
\put(10.4,4.6){$h = \frac{1}{200}$}
\end{picture}
\caption{ \small \label{Fig:A2Dres10}
Test (\ref{eq:Adv2D}) with flux-based decomposition.
Maximum errors for the schemes
TW2 ({\tiny $\square$} marks), CS2 ($\circ$ marks) and SH2 ($\diamond$ marks)
as function of Courant numbers $\nu =  2\pi\dt/h$, between $\frac{1}{2}$ 
and $2$,
for the grids with $h = \frac{1}{50}$ (left), $h = \frac{1}{100}$ (middle)
and $h = \frac{1}{200}$ (right).
}
\vspace{-2mm}
\end{center}
\end{figure}

The errors in the maximum norm are given in
Figure~\ref{Fig:A2Dres10} on three uniform grids, $\dx = \dy = h$ with
$h = \frac{1}{50}, \frac{1}{100}, \frac{1}{200}$, 
and this is to be compared with the 
results in Figure~\ref{Fig:A2Dres00}. Again the dashed line 
in the figure gives the result for the scheme where in each time 
step the explicit trapezoidal rule is applied twice, with step-size 
$\frac{1}{2}\dt$, over the whole region.

Compared to Figure~\ref{Fig:A2Dres00}, the negative effect of the flux-based
splitting on the accuracy of the CS2 scheme is obvious. Here it is to be
noted that the vertical axis in Figure~\ref{Fig:A2Dres10} has been shifted to
include the error lines in the plots.

However, the accuracy of the TW2 and SH2 schemes has deteriorated as well, 
which is most clear by comparing these results with the ones for the explicit 
trapezoidal rule with small step-size $\frac{1}{2}\dt$ over the whole region,
which we may consider here as a 'target' solution.
In contrast to Figure~\ref{Fig:A2Dres00}, where the results of the
TW2 and SH2 schemes were close to these reference solutions, now the
errors with the TW2 and SH2 schemes are much larger, in particular on the 
finer grids.

Instead of errors $C \dt^2$ with a fixed constant $C$, the constants in 
front of the global errors are now proportional to $h^{-1}$, and 
comparing the results on the three grids for fixed ratios $\dt/h$, it can
be observed that the order of convergence for TW2 and SH2 is now only one.
So we still have convergence with these schemes in the maximum
norm, but it is much slower than for the cell-based splittings.

The largest errors are found near the interfaces. Measuring the errors
in the $L_1$-norm would yield convergence with one order higher, similar
as for the 1D test in Table~\ref{Tab:Test1FluxBased}.
Convergence with order two in the $L_1$-norm with the TW2 and SH2 schemes 
may be satisfactory for many applications.

\section{Conclusions}

In this paper the accuracy of partitioned Runge-Kutta methods has been 
studied, with applications to explicit multirate schemes. When such 
methods are applied to PDEs, it is not sufficient to look at the order 
for non-stiff problems. The interfaces between the regions where different 
methods --\,or different time steps\,-- are applied act like
time dependent boundary conditions, and order reduction is to be
expected.

To see the effect of the partitioning at the interfaces, the accuracy of
the schemes was mainly considered in the maximum norm. Convergence 
in the discrete $L_1$-norm is in general one order larger. This due to
fact that the largest errors are confined to small spatial regions near
the interfaces.

To guarantee mass conservation during all stages of the computation,
a decomposition based on fluxes seems attractive. However, it was seen that
the order of convergence for smooth problems will be smaller compared to
cell-based splittings.
On the other hand, for partitioned Runge-Kutta methods with different 
weights, the cell-based splittings may lead to an incorrect propagation
of discontinuities.

If a high accuracy is required, then one would like to use high-order
methods, of course, and the decompositions considered in this paper do
not seem to be very suited.
Alternatives are the use of smooth partitions of unity, similar to the 
approach in \cite{MPRW98, PBJ04b} for parabolic problems, or an approach 
with overlapping regions.
The study of convergence and monotonicity/SSP properties of such methods
is part of our current research.

\bigskip\noindent
{\bf Acknowledgement.}
This paper originated from work of W.\,H.\ with Anna Mozartova and 
Valeriu Savcenco. The contributions of A.M.\ and V.S.\, on the design 
of multirate methods and monotonicity properties of these methods, are 
contained in \cite{HMS13}. 
They are thanked for helpful comments on preliminary convergence results 
for first-order methods.

\end{document}